\documentclass[12pt, A4paper]{amsart}
  \usepackage{latexsym}
  \usepackage[all]{xy}
  \usepackage{amsthm}
 \numberwithin{equation}{section}
  \usepackage{color}
  \usepackage{hyperref}

\newcommand{\Ra}{\Rightarrow}
\def\cK{\mathcal K}
\def\ot{{\otimes}}
\def\1c#1{\stackrel{#1}{\to}}
\def\2c#1{\stackrel{#1}{\Ra}}

\def\Mnd{\mathrm{Mnd}}
\def\Cmd{\mathrm{Cmd}}
\def\Entw{\mathrm{Entw}}

  \newtheorem{proposition}{Proposition}[section]

  \newtheorem{corollary}[proposition]{Corollary}
  \newtheorem{theorem}[proposition]{Theorem}

  \theoremstyle{definition}
  
  \newtheorem{example}[proposition]{Example}

  \theoremstyle{remark}

  \newcounter{c}
  
  \newcommand{\etyk}[1]{\vspace{-7.4mm}$$\begin{equation}\Label{#1}
  \addtocounter{c}{1}}
  \renewcommand{\]}{\ifnum \value{c}=1 $$\else \end{equation}\fi}
\begin{document}

 \title{The 2-category of weak entwining structures}
 \author{Gabriella B\"ohm}
 \address{Research Institute for Particle and Nuclear Physics, Budapest,
 \newline\indent H-1525
 Budapest 114, P.O.B.\ 49, Hungary}
  \email{G.Bohm@rmki.kfki.hu}
  \begin{abstract}
A weak entwining structure in a 2-category $\cK$ consists of a monad $t$ and a
comonad $c$, together with a 2-cell relating both structures in a way which
generalizes a mixed distributive law. 
A weak entwining structure can be characterized as a compatible pair of a
monad and a comonad, in 2-categories generalizing the 2-category of comonads
and the 2-category of monads in $\cK$, respectively. This observation is used
to define a 2-category $\Entw^w(\cK)$ of weak entwining structures in $\cK$.  
If the 2-category $\cK$ admits Eilenberg-Moore constructions for both monads
and comonads and idempotent 2-cells in $\cK$ split, then there are
pseudo-functors from $\Entw^w(\cK)$ to the 2-category of monads and to the
2-category of comonads in $\cK$, taking a weak entwining structure $(t,c)$ to
a `weak lifting' of $t$ for $c$ and a `weak lifting' of $c$ for $t$,
respectively. The Eilenberg-Moore objects of the lifted monad and the lifted
comonad are shown to be equivalent. 
If $\cK$ is the 2-category of functors induced by bimodules, then these 
Eilenberg-Moore objects are isomorphic to the usual category of weak entwined
modules. 
\end{abstract}
  \maketitle

\section*{Introduction}

Mixed distributive laws \cite{Beck} in a 2-category $\cK$ (or `entwining
structures', as they are called more often in the Hopf algebraic terminology),
can be described in some equivalent ways \cite{Street}. They are monads in the 
2-category $\Cmd(\cK)$ of comonads in $\cK$, equivalently, they are comonads
in the 2-category $\Mnd(\cK)$ of monads in $\cK$. Consequently, they can be
regarded as 0-cells of a 2-category $\Entw(\cK)$, defined to be isomorphic to
$\Mnd(\Cmd(\cK))\cong \Cmd(\Mnd(\cK))$.  

If a 2-category $\cK$ admits Eilenberg-Moore constructions for monads,
that is, the inclusion 2-functor $I:\cK \to \Mnd(\cK)$ possesses a right
2-adjoint $J$, then the 2-functor $\Cmd(J)$ takes a mixed distributive law of
a monad $t$ and a comonad $c$ in $\cK$ to a comonad $J(t) \1c {\overline c}
J(t)$, which is a lifting of $c$, cf. \cite{PowWat}.
Symmetrically, if $\cK$ admits Eilenberg-Moore constructions 
for comonads, that is, the inclusion 2-functor $I_*:\cK \to \Cmd(\cK)$ 
possesses a right 2-adjoint $J_*$, then $\Mnd(J_*)$ takes $(t,c)$ to a monad
$J_*(c) \1c {\overline t} J_*(c)$, which is a lifting of $t$.
If Eilenberg-Moore constructions in $\cK$ exist both for monads and comonads,
then the 2-functors $J_* \Cmd(J)$ and $J\Mnd(J_*)$ are 2-naturally
isomorphic. In particular, the lifted monad ${\overline t}$ and the lifted
comonad ${\overline c}$ possess isomorphic Eilenberg-Moore objects, see
\cite{PowWat}. In the case when $\cK$ is the 2-category $\mathrm{CAT}=${\tt
[Categories; Functors; Natural Transformations]}, this is the category of 
$(t,c)$-bimodules, also called `entwined modules'. 

In order to treat algebra extensions by weak bialgebras in \cite{BNSz},
entwining structures were generalized to `weak entwining structures' in
\cite{CaeDeG}. A weak entwining structure in a 2-category $\cK$ also consists
of a monad $t$ and a comonad $c$, together with a 2-cell $tc \2c {} ct$, but
the compatibility axioms with the unit of the monad and the counit of the
comonad are weakened.
We are not aware of any characterization of a weak entwining structure as a
monad or as a comonad in some 2-category. Instead, in this note we observe that
a weak entwining structure in an arbitrary 2-category $\cK$ can be described
as a compatible pair of a comonad in a 2-category $\Mnd^\iota(\cK)$, which
extends $\Mnd(\cK)$, and a monad in $\Cmd^\pi(\cK):=\Mnd^\iota(\cK_*)_*$
(where $(-)_*$ means the vertically opposite 2-category). 
This observation is used to define in Section \ref{sec:wentw} a 2-category
$\Entw^w(\cK)$, whose 0-cells are weak entwining structures in $\cK$ and
whose 1-cells and 2-cells are also compatible pairs of 1-cells and 2-cells,
respectively, in $\Mnd(\Cmd^\pi(\cK))$ and $\Cmd(\Mnd^\iota(\cK))$. By 
construction, the 2-category $\Entw^w(\cK)$ comes equipped with 2-functors
$A:\Entw^w(\cK) \to \Cmd(\Mnd^\iota(\cK))$ and $B:\Entw^w(\cK) \to
\Mnd(\Cmd^\pi(\cK))$.   

If a 2-category $\cK$ admits Eilenberg-Moore constructions for monads and
idempotent 2-cells in $\cK$ split, then the 2-functor $J$ above factorizes
through the inclusion $\Mnd(\cK) \hookrightarrow \Mnd^\iota(\cK)$ and an
appropriate pseudo-functor $Q:\Mnd^\iota(\cK) \to \cK$. The image of a weak
entwining structure $(t,c)$ under the pseudo-functor $\Cmd(Q)A$ is a `weak
lifting' of $c$ for $t$, cf. \cite{weak_EM}. Symmetrically, if $\cK$ admits
Eilenberg-Moore constructions for comonads and idempotent 2-cells in $\cK$
split, then there is a pseudo-functor $Q_*:\Cmd^\pi(\cK) \to \cK$, such that
$\Mnd(Q_*)B$ takes a weak entwining structure $(t,c)$ to a weak lifting of $t$
for $c$. If Eilenberg-Moore constructions in $\cK$ exist both for monads and
comonads and also idempotent 2-cells in $\cK$ split, then we prove in Section
\ref{sec:EM_iso} that the pseudo-functors $J_* \Cmd(Q)A$ and $J \Mnd(Q_*)B:
\Entw^w(\cK)\to \cK$ are pseudo-naturally equivalent. In particular, for any
weak entwining structure $(t,c)$, the weak lifting of $t$ for $c$, and the
weak lifting of $c$ for $t$, possess equivalent Eilenberg-Moore objects.  

As a motivating example, we can consider the 2-category $\cK$ obtained as the
image of the bicategory $\mathrm{BIM}_k=${\tt [Algebras; Bimodules; Bimodule
Maps]} (over a commutative ring $k$) under the hom 2-functor
$\mathrm{BIM}_k(k,-):\mathrm{BIM}_k \to \mathrm{CAT}$. A weak entwining
structure $((-)\ot_R T, (-)\ot_R C)$ in this 2-category is given by a
$k$-algebra $R$, an $R$-ring $T$, an $R$-coring $C$ and an $R$-bimodule map $C
\otimes_R T \to T \otimes_R C$. 
In this case, we obtain that the Eilenberg-Moore category of the weakly lifted 
comonad $\overline{(-) \otimes_R C}$ (on the category $M_T$ of $T$-modules)
is isomorphic to the Eilenberg-Moore category of the weakly lifted monad
$\overline{(-) \otimes_R T}$ (on the category $M^C$ of $C$-comodules), and it
is  isomorphic also to  $\Entw^w(\cK)((M_k,M_k),((-)\ot_R T,(-)\ot_R C))$,
known as the category of `weak entwined modules'. In particular, if $R$ is a
trivial $k$-algebra (i.e. $R=k$), we re-obtain \cite[Proposition
  2.3]{Brz:coring}.  
\bigskip

{\bf Notations.}
We assume that the reader is familiar with the theory of 2-categories. For a
review of the occurring notions (such as a 2-category, a 2-functor and a
2-adjunction, monads, adjunctions and Eilenberg-Moore construction in a
2-category) we refer to the article \cite{KeSt}.

In a 2-category $\cK$, horizontal composition is denoted by juxtaposition and
vertical composition is denoted by $\ast$, 1-cells are represented by an arrow
$\1c{}$ and 2-cells are represented by $\2c {}$. 

For any 2-category $\cK$, $\Mnd(\cK)$ denotes the 2-category of monads in
$\cK$ as in \cite{Street} and $\Cmd(\cK):=\Mnd(\cK_*)_*$ denotes the
2-category of comonads in $\cK$, where $(-)_*$ refers to the vertical opposite
of a 2-category. Throughout, we denote by $I:\cK\to \Mnd(\cK)$ the inclusion
2-functor (with underlying maps $k\mapsto (k,k,k)$, $V\mapsto (V,V)$, $\omega
\mapsto \omega$ on the 0-, 1-, and 2-cells, respectively). Its right
2-adjoint, if it exists, is denoted by $J$. The inclusion 2-functor $\cK
\to \Cmd(\cK)$ is denoted by $I_*$ and its right 2-adjoint, whenever it
exists, is denoted by $J_*$. 

If a 2-category $\cK$ admits Eilenberg-Moore constructions for monads
(i.e. the 2-functor $J$ exists), then 
any monad $(k\1c t k, tt\2c \mu t,k\2c \eta t)$ in $\cK$ determines a
canonical adjunction $(k\1c f J(t), J(t)\1c v k, fv \2c \epsilon J(t), k \2c
\eta vf)$ such that $(t,\mu,\eta)=(vf,v\epsilon f,\eta)$, cf. \cite[Theorem
2]{Street}. Throughout, these notations are used for this canonical
adjunction. For a monad $(t',\mu',\eta')$, the 
canonical adjunction is denoted by $(f',v',\epsilon',\eta')$, etc.

We say that in a 2-category $\cK$ idempotent 2-cells split if, for any
2-cell $V \2c \Theta V$ in $\cK$ such that $\Theta\ast \Theta = \Theta$, there
exist a 1-cell ${\widehat V}$ and 2-cells $V \2c \pi {\widehat V}$ and
${\widehat V} \2c \iota V$, such that $\pi \ast \iota ={\widehat V}$ and
$\iota \ast \pi =\Theta$.

\section{The 2-category of weak entwining structures}
\label{sec:wentw}

Consider a monad $(k\1c t k,tt \2c \mu t,k \2c \eta t)$ and a comonad $(k\1c c
k, c\2c \delta cc, c\2c \varepsilon k)$ in a 2-category $\cK$ and a 2-cell $tc
\2c \psi ct$. 
The triple $(t,c,\psi)$ is termed a {\em weak entwining structure} provided
that the following axioms in \cite{CaeDeG} hold.
\begin{eqnarray}
&& \psi \ast \mu c 
= c\mu \ast \psi t \ast t\psi;
\label{eq:wentw_m}\\
&& \delta t \ast \psi 
= c\psi \ast \psi c \ast t\delta;
\label{eq:wentw_d}\\
&& \psi \ast \eta c 
= c\varepsilon t \ast c\psi \ast c\eta c \ast \delta;
\label{eq:wentw_u}\\
&& \varepsilon t \ast \psi 
= \mu \ast t\varepsilon t \ast t \psi \ast t\eta c.
\label{eq:wentw_e}
\end{eqnarray}
The most important difference between such a weak entwining structure and a
usual entwining structure (i.e. mixed distributive law) is that in the weak
case $(c,\psi)$ is no longer a 1-cell $t \1c{} t$ in $\Mnd(\cK)$ and
$(t,\psi)$ is not a 1-cell $c \1c{} c$ in $\Cmd(\cK)$. Still, as it was
observed in \cite{weak_EM}, $(t\1c{(c,\psi)} t,\mu,\eta)$ is a monad and $(c
\1c{(t,\psi)} c,\delta,\varepsilon)$ is a comonad in an extended 2-category of
(co)monads in $\cK$, recalled in the following theorem.

\begin{theorem}[\cite{weak_EM}, Corollary 1.4 and Theorem 3.5]
\label{thm:Mnd^i} 
For any 2-category $\cK$, the following data constitute a 2-category, to be
denoted by $\Mnd^\iota(\cK)$.
\begin{itemize}
\item[] \underline{0-cells} are monads $(k \1c t k,\mu,\eta)$ in $\cK$.
\item[] \underline{1-cells} $(k\1c t k,\mu,\eta)\1c{(V,\psi)} (k' \1c {t'}
  k',\mu',\eta')$ are pairs, consisting of a 1-cell $k\1c V k'$ and a 2-cell
  $t'V \2c \psi Vt$ in $\cK$ such that 
\begin{equation}\label{eq:1-cell}
V\mu \ast \psi t \ast t' \psi = \psi \ast \mu' V.
\end{equation}
\item[] \underline{2-cells} $(V,\psi)\2c{\omega} (W,\phi)$ are 2-cells
$V\2c{\omega}W$ in $\cK$, satisfying 
\begin{equation}\label{eq:2-cell}
\omega t\ast \psi = W\mu \ast \phi t \ast t'\omega t \ast t' \psi
  \ast t' \eta' V.
\end{equation}
\item[] \underline{Horizontal and vertical compositions} are the same as in
  $\cK$. 
\end{itemize}
The 2-category $\Mnd^\iota(\cK)$ contains $\Mnd(\cK)$ as a vertically full
2-subcategory. 

Moreover, if $\cK$ admits Eilenberg-Moore constructions for monads and
idempotent 2-cells in $\cK$ split, then the following maps determine a
 pseudo-functor 
$Q:\Mnd^\iota(\cK)\to \cK$.
\begin{itemize}
\item[] For a \underline{0-cell} $(t,\mu,\eta)$,
  $Q(t,\mu,\eta):=J(t,\mu,\eta)$. 
\item[] For a \underline{1-cell} $(t,\mu,\eta)\1c{(V,\psi)} (t',\mu',\eta')$,
$Q(V,\psi)$ is the unique 1-cell $Q(t,\mu,\eta)\1c{} Q(t',\mu',$ $\eta')$ in
  $\cK$ for which 
\begin{equation}\label{eq:Q_1-cell}
v'\epsilon' Q(V,\psi)=\pi \ast Vv \epsilon \ast \psi v \ast t'\iota.
\end{equation}
\item[]For a \underline{2-cell} $(V,\psi)\2c{\omega} (W,\phi)$, $Q(\omega)$
  is the unique 2-cell $Q(V,\psi)\2c{} Q(W,\phi)$ in $\cK$ for which
\begin{equation}\label{eq:Q_2-cell}
v'Q(\omega)= \pi \ast \omega v \ast \iota,
\end{equation}
\end{itemize}
 where $Vv \2c \pi v'Q(V,\psi) \2c \iota Vv$ denote a chosen splitting 
of the idempotent 2-cell 
\begin{equation}\label{eq:idemp}
Vv\epsilon \ast \psi v\ast \eta'Vv:Vv\2c{}Vv,
\end{equation} 
for any 1-cell $(V,\psi)$ in $\Mnd^\iota(\cK)$.

For 1-cells $t \1c {(V,\psi)} t' \1c {(V',\psi')} t''$ in $\Mnd^\iota(\cK)$,
the coherence natural iso 2-cell \break
$Q((V',\psi')(V,\psi)) \2c \cong Q(V',\psi')Q(V,\psi)$ is the unique 2-cell
$\gamma$ for that  \break
$v''\gamma = 
\big( v''Q((V',\psi')(V,\psi)) \2c \iota V'Vv
\2c {V' \pi} V'v' 
Q(V,\psi) 
\2c {\pi Q(V,\psi)} 
v'' Q(V',\psi')Q(V,\psi) \big)$
(so
$v''\gamma^{-1} = 
\big( v'' Q(V',\psi')Q(V,\psi) 
\2c {\iota Q(V,\psi)} 
V'v' Q(V,\psi) 
\2c {V' \iota } 
V'Vv \2c \pi v''Q((V',\psi')(V,\psi)) \big)$).
With the convention of choosing a trivial splitting $Vv \2c {Vv} Vv \2c {Vv}
Vv$ whenever \eqref{eq:idemp} is an identity 2-cell, the image of any identity
1-cell $t \1c {(k,t)} t$ under $Q$ becomes equal to the identity 1-cell
$Q(t)$. This convention also ensures that the composite pseudo-functor
$\Mnd(\cK)\hookrightarrow \Mnd^\iota(\cK)\1c Q \cK$ is equal to $J$.
The pseudo-natural isomorphism class of $Q$ does not depend on the choice of
the 2-cells $\pi$ and $\iota$. 
\end{theorem}

For any 2-category $\cK$, we put $\Cmd^\pi(\cK):=\Mnd^\iota(\cK_*)_*$. Applying
Theorem \ref{thm:Mnd^i} to the 2-category $\cK_*$, we conclude that whenever
$\cK$ admits Eilenberg-Moore constructions for comonads and idempotent 2-cells
in $\cK$ split, $J_*$ extends to a pseudo-functor $Q_*:\Cmd^\pi(\cK)\to \cK$.

After all these preparations, we are ready to construct a 2-category of weak
entwining structures in any 2-category $\cK$.

\begin{theorem}\label{thm:wentw}
For any 2-category $\cK$, the following data constitute a 2-category, to be
denoted by $\Entw^w(\cK)$.
\begin{itemize}
\item[] \underline{0-cells} are triples $((k \1c t k,\mu,\eta),(k\1c c
  k,\delta,\varepsilon),\psi)$, consisting of a monad $(k \1c t k,\mu,\eta)$,
  a comonad $(k\1c c k,\delta,\varepsilon) $ and a 2-cell $tc\2c \psi ct$ in
  $\cK$, such that   
\begin{itemize}
\item[$\bullet$]$(t\1c{(c,\psi)}t,\delta,\varepsilon)$ is a comonad in
  $\Mnd^\iota(\cK)$ and 
\item[$\bullet$] $(c\1c{(t,\psi)}c,\mu,\eta)$ is a monad in $\Cmd^\pi(\cK)$. 
\end{itemize}
\item[] \underline{1-cells} 
$((k \1c t k,\mu,\eta),(k\1c c k,\delta,\varepsilon),\psi)
  \1c{(W,\alpha,\beta)} ((k' \1c {t'} k',\mu',\eta'),(k'\1c {c'}
  k',\delta',\varepsilon'),\psi')$ are triples, consisting of a 
  1-cell $k\1c W k'$ and 2-cells $t'W \2c \alpha Wt$ and $Wc\2c \beta c'W$ in
  $\cK$, such that 
\begin{itemize}
\item[$\bullet$] $(t\1c{(c,\psi)}t,\delta,\varepsilon) \1c{((W,\alpha),\beta)}
  (t'\1c{(c',\psi')}t',\delta',\varepsilon')$ is a 1-cell in
  $\Cmd(\Mnd^\iota(\cK))$ and  
\item[$\bullet$] $(c\1c{(t,\psi)}c,\mu,\eta) \1c{((W,\beta),\alpha)}
  (c'\1c{(t',\psi')}c',\mu',\eta')$ is a 1-cell in $\Mnd(\Cmd^\pi(\cK))$.
\end{itemize}
\item[] \underline{2-cells} $(W,\alpha,\beta)\2c{\omega} (W',\alpha',\beta')$
  are 2-cells $W\2c{\omega}W'$ in $\cK$, such that
\begin{itemize}
\item[$\bullet$] $((W,\alpha),\beta)\2c \omega ((W',\alpha'),\beta')$ is a
  2-cell in $\Cmd(\Mnd^\iota(\cK))$ and  
\item[$\bullet$] $((W,\beta),\alpha)\2c \omega ((W',\beta'),\alpha')$ is a
  2-cell in $\Mnd(\Cmd^\pi(\cK))$.  
\end{itemize}
\item[] \underline{Horizontal and vertical compositions} are the same as in
  $\cK$. 
\end{itemize}
\end{theorem}

\begin{proof}
In order to see that 0-cells in $\Entw^w(\cK)$ are precisely the weak
entwining structures, 
note that \eqref{eq:wentw_m} expresses the requirement that $t\1c{(c,\psi)}t$
is a 1-cell in $\Mnd^\iota(\cK)$ and \eqref{eq:wentw_d} means that $c
\1c{(t,\psi)} c$ is a 1-cell in $\Cmd^\pi(\cK)$. Axiom \eqref{eq:wentw_u} 
means that $(k,c)\2c \eta (t,\psi)$ is a 2-cell in $\Cmd^\pi(\cK)$ and 
\eqref{eq:wentw_e} holds if and only if $(c,\psi)\2c \varepsilon (k,t)$ is a
2-cell in $\Mnd^\iota(\cK)$.   
If these four conditions hold, then also $(t,\psi)(t,\psi)\2c \mu
(t,\psi)$ is a 2-cell in $\Cmd^\pi(\cK)$. That is,
\begin{eqnarray*}
c\varepsilon t \ast c\psi \ast c\mu c \ast \psi tc \ast t \psi c \ast tt\delta 
&\stackrel{\eqref{eq:wentw_m}}{=}&
c \varepsilon t\ast c\psi \ast \psi c \ast \mu cc \ast tt \delta 
=c \varepsilon t\ast c\psi \ast \psi c \ast t \delta \ast \mu c\\
&\stackrel{\eqref{eq:wentw_d}}{=}&
c \varepsilon t \ast  \delta t \ast \psi \ast  \mu c
= \psi \ast  \mu c.
\end{eqnarray*}
Similarly, (\ref{eq:wentw_m}-\ref{eq:wentw_e}) imply that $(c,\psi)\2c \delta
(c,\psi)(c,\psi)$ is a 2-cell in $\Mnd^\iota(\cK)$, i.e. 
\begin{eqnarray*}
cc \mu \ast c\psi t \ast \psi ct \ast t \delta t \ast t \psi \ast t \eta c
&\stackrel{\eqref{eq:wentw_d}}{=}&
cc \mu  \ast  \delta tt \ast \psi t \ast t\psi \ast t \eta c
= \delta t \ast c \mu  \ast \psi t \ast t\psi \ast t \eta c\\
&\stackrel{\eqref{eq:wentw_m}}{=}&
 \delta t \ast \psi \ast \mu c \ast t \eta c
=\delta t \ast \psi.
\end{eqnarray*}

By Theorem \ref{thm:Mnd^i}, a triple $(k \1c{W}k', t'W \2c \alpha Wt,Wc \2c
\beta c'W)$ is a 1-cell
$((k \1c t k,\mu ,\eta ),(k\1c c k,\delta,\varepsilon),\psi) \1c{} 
((k' \1c {t'} k',\mu',\eta'),(k'\1c {c'} k',\delta',\varepsilon'),\psi')$ in 
$\Entw^w(\cK)$ if and only if the following equalities hold. 
\begin{eqnarray}
&&\alpha \ast \mu'W=W \mu \ast \alpha t \ast t' \alpha;
\label{eq:1-cell_7}\\
&&\alpha \ast \eta' W = W\eta;
\label{eq:1-cell_12}\\
&&\delta'W \ast \beta = c' \beta \ast \beta c \ast W\delta;
\label{eq:1-cell_9}\\
&&\varepsilon' W \ast \beta = W\varepsilon;
\label{eq:1-cell_10}\\
&&c'W\mu \ast c'\alpha t \ast \psi' Wt \ast t'\beta t \ast t'W\psi \ast t'
W\eta c = \beta t \ast W \psi \ast \alpha c
\label{eq:1-cell_8}\\
&&c'W \varepsilon t \ast c' W\psi \ast c'\alpha c \ast \psi'Wc \ast t'\beta c
\ast t' W \delta = \beta t \ast W \psi \ast \alpha c.
\label{eq:1-cell_11}
\end{eqnarray}
The equality \eqref{eq:1-cell_7} is equivalent to saying that $t
\1c{(W,\alpha)} t'$ is a 1-cell in $\Mnd^\iota(\cK)$ and \eqref{eq:1-cell_9} is
equivalent to $c \1c{(W,\beta)} c'$ being a 1-cell in $\Cmd^\pi(\cK)$.
The equality \eqref{eq:1-cell_11} means (after being simplified using
\eqref{eq:1-cell_10}) that $(t',\psi')(W,\beta) \2c \alpha
(W,\beta)(t,\psi)$ is a 2-cell in $\Cmd^\pi(\cK)$ and \eqref{eq:1-cell_8} means
(after being simplified using \eqref{eq:1-cell_12}) 
that $(W,\alpha)(c,\psi)\2c \beta (c',\psi')(W,\alpha)$ is a 2-cell in
$\Mnd^\iota(\cK)$. 
Conditions \eqref{eq:1-cell_7} and \eqref{eq:1-cell_12} mean that
$(c \1c{(t,\psi)}c,\mu,\eta) \2c{((W,\beta),\alpha)} (c'
\1c{(t',\psi')}c',\mu',\eta')$ is a 2-cell in $\Mnd(\Cmd^\pi(\cK))$, while  
\eqref{eq:1-cell_9} and \eqref{eq:1-cell_10} express that 
$(t\1c{(c,\psi)}t,\delta,\varepsilon)\2c{((W,\alpha),\beta)}
(t'\1c{(c',\psi')}t',\delta',\varepsilon')$ is a 2-cell in
$\Cmd(\Mnd^\iota(\cK))$.  

A 2-cell $W \2c \omega W'$ in $\cK$ is a 2-cell $(W,\alpha,\beta) \2c {}
(W',\alpha',\beta')$ in $\Entw^w(\cK)$ if and only if 
\begin{eqnarray}
&&\alpha' \ast t'\omega = \omega t \ast \alpha
\label{eq:2-cell_14}\\
&&\beta' \ast \omega c=c'\omega \ast \beta .
\label{eq:2-cell_13}
\end{eqnarray}

For any weak entwining structure $((k\1c t k,\mu,\eta),(k\1c c
k,\delta,\varepsilon),\psi)$ in $\cK$, the triple $(W=k,\alpha=t,\beta=c)$
satisfies the equalities (\ref{eq:1-cell_7}-\ref{eq:1-cell_11}). Hence it is
an (identity) 1-cell in $\Entw^w(\cK)$.  
The sets of 1-cells and 2-cells in $\Cmd(\Mnd^\iota(\cK))$ and
$\Mnd(\Cmd^\pi(\cK))$ are closed under the horizontal composition in $\cK$ by
Theorem \ref{thm:Mnd^i}. Therefore the horizontal composite of 1-cells and
2-cells in $\Entw^w(\cK)$ is a 1-cell and a 2-cell in $\Entw^w(\cK)$,
respectively. 

For any 1-cell $(W,\alpha,\beta)$ in $\Entw^w(\cK)$, the identity 2-cell $W\2c
W W$ in $\cK$ satisfies \eqref{eq:2-cell_14} and \eqref{eq:2-cell_13}. Hence
it is an (identity) 2-cell in $\Entw^w(\cK)$. Since the sets of 2-cells in
$\Cmd(\Mnd^\iota(\cK))$ and $\Mnd(\Cmd^\pi(\cK))$ are closed under the 
vertical composition in $\cK$ by Theorem \ref{thm:Mnd^i}, the vertical
composite of 2-cells in $\Entw^w(\cK)$ is a 2-cell in $\Entw^w(\cK)$
again. 

Associativity and unitality of the horizontal and vertical compositions in
$\Entw^w(\cK)$ and the interchange law follow by the respective properties of
$\cK$. 
\end{proof}

From Theorem \ref{thm:wentw}, we immediately deduce the existence of some
2-functors. 

\begin{corollary}\label{cor:functors}
For any 2-category $\cK$, the following assertions hold.
\begin{itemize}
\item[{(1)}] There is a 2-functor $Y:\cK\to \Entw^w(\cK)$, determined by the
  maps $k\mapsto (I(k),I_*(k),k)$, $V\to (V,V,V)$ and $\omega \mapsto \omega$
  on the 0-, 1-, and 2-cells, respectively.
\item[{(2)}] There is a 2-category isomorphism $\Phi:\Entw^w(\cK)\cong
  \Entw^w(\cK_*)_*$, determined by the maps $(t,c,\psi)\mapsto (c,t,\psi)$,
  $(W,\alpha,\beta)\mapsto (W,\beta,\alpha)$ and $\omega\mapsto \omega$ on the
  0-, 1-, and 2-cells, respectively. In particular, for any weak entwining
  structures $(t,c,\psi)$ and $(t',c',\psi')$ in $\cK$, there is a category
  isomorphism  $\Entw^w(\cK)((t,c,\psi),(t',c',\psi'))\cong$ 
$\Entw^w(\cK_*)_*((c,t,\psi),(c',t',\psi'))$, which is 2-natural both in
$(t,c,\psi)$ and $(t',c',\psi')$. 
\item[{(3)}] There is a 2-functor $A:\Entw^w(\cK) \to \Cmd(\Mnd^\iota(\cK))$,
  determined by the maps $((t,\mu,\eta),(c,\delta,\varepsilon),\psi)\mapsto
  (t\1c{(c,\psi)}t,\delta,\varepsilon)$, $(W,\alpha,\beta)\mapsto
  ((W,\alpha),\beta)$ and 
  $\omega \mapsto \omega$ on the 0-, 1-, and 2-cells, respectively. 
\item[{(4)}] There is a 2-functor $B:\Entw^w(\cK) \to \Mnd(\Cmd^\pi(\cK))$,
  determined by the maps $((t,\mu,\eta),(c,\delta,\varepsilon),\psi)\mapsto
  (c\1c{(t,\psi)}c,\mu,\eta)$, $(W,\alpha,\beta)\mapsto ((W,\beta),\alpha)$
  and $\omega \mapsto \omega$ on the 0-, 1-, and 2-cells, respectively. 
\end{itemize}
\end{corollary}

In contrast to the case of usual entwining structures, there seems to be no
reason to expect that the 2-functors $A$ and $B$ in Corollary
\ref{cor:functors} are isomorphisms.

\section{ Equivalence  of Eilenberg-Moore objects}
\label{sec:EM_iso}

If a 2-category $\cK$ admits Eilenberg-Moore constructions for both monads and
comonads and idempotent 2-cells in $\cK$ split, then by Theorem
\ref{thm:Mnd^i} and Corollary \ref{cor:functors}, there are two
pseudo-functors $J_* \Cmd(Q)A$ and $J\Mnd(Q_*)B:\Entw^w(\cK)\to \cK$. The aim
of this section is to prove that both are right biadjoints of $Y$ in Corollary
\ref{cor:functors}(1), hence they are pseudo-naturally equivalent. 
Consequently, for any weak entwining structure $(t,c,\psi)$ in $\cK$, the monad
$Q_*(c\1c{(t,\psi)}c)$ and the comonad $Q(t\1c{(c,\psi)}t)$ in $\cK$ possess
equivalent Eilenberg-Moore objects. 

Recall that any pseudo-functor $Q:{\mathcal A} \to {\mathcal B}$ between
2-categories, induces a pseudo-functor $\mathrm{Cmd}(Q):\mathrm{Cmd}({\mathcal
A}) \to \mathrm{Cmd}({\mathcal B})$ with underlying maps as follows. A comonad
$(A \1c c A, \delta, \varepsilon )$ in ${\mathcal A}$ is taken to the comonad
$Q(A)\1c {Q(c)} Q(A)$, with comultiplication $Q(c) \2c {Q(\delta)} Q(cc) \2c
\cong Q(c)Q(c)$ and counit $Q(c) \2c {Q(\varepsilon)} Q(1_A) \2c \cong
1_{Q(A)}$. A 1-cell $(A \1c c A, \delta, \varepsilon ) \1c {(V,\psi)} (A' \1c
{c'} A', \delta', \varepsilon')$ in $\mathrm{Cmd}({\mathcal A})$ is taken to a
pair consisting of the 1-cell $Q(A)\1c {Q(V)} Q(A')$ and the 2-cell $Q(V)Q(c)
\2c \cong Q(Vc) \2c {Q(\psi)} Q(c'V) \2c \cong Q(c') Q(V)$ in ${\mathcal
B}$. A 2-cell $\omega$ in $\mathrm{Cmd}({\mathcal A})$ is taken to
$Q(\omega)$ . $\mathrm{Cmd}(Q)$ is a pseudo-functor with the same coherence
isomorphisms as $Q$.

\begin{proposition}\label{prop:comodules}
Consider a 2-category $\cK$ which admits Eilenberg-Moore constructions for
monads and in which idempotent 2-cells split. 
Let $l$ be a 0-cell and $((k \1c t k,\mu,\eta),(k\1c c
k,\delta,\varepsilon),\psi)$ be weak entwining structure in $\cK$.
The following categories are isomorphic.
\begin{itemize}
\item[{(i)}] The Eilenberg-Moore category
$\Cmd(\cK)(I_*(l),\Cmd(Q)(t\1c{(c,\psi)}t,\delta,\varepsilon))$ of the
  co\-monad $\cK(l,Q(t\1c{(c,\psi)}t)):\cK(l,Q(t))\to \cK(l,Q(t))$;
\item[{(ii)}] the category $\Entw^w(\cK)(Y(l),(t,c,\psi))$.
\end{itemize}
Moreover, these isomorphisms provide the 1-cell parts of a pseudo-natural
isomorphism $\Cmd(\cK)(I_*(-),\Cmd(Q)A(-)) \cong \Entw^w(\cK)(Y(-),-)$. 
\end{proposition}

\begin{proof}
By (\ref{eq:1-cell_7}-\ref{eq:1-cell_11}), the objects in the category
$\Entw^w(\cK)(Y(l),(t,c,\psi))$ are triples $(l\1c W k,tW\2c \varrho W,W\2c
\kappa cW)$, such that $I(l) \1c{(W,\varrho)} t$ is a 1-cell in $\Mnd(\cK)$,
$I_*(l) \1c{(W,\kappa)} c$ is a 1-cell in $\Cmd(\cK)$ and 
\begin{equation}\label{eq:entw_mod}
c\varrho \ast \psi W \ast t\kappa = \kappa \ast \varrho.
\end{equation}
Morphisms $(W,\varrho,\kappa)\to (W',\varrho',\kappa')$ in
$\Entw^w(\cK)(Y(l),(t,c,\psi))$ are 2-cells $W\2c \omega W'$ in $\cK$, such
that $(W,\varrho) \2c{\omega} (W',\varrho')$ is a 2-cell in $\Mnd(\cK)$ and
$(W,\kappa)\2c \omega (W',\kappa')$ is a 2-cell in $\Cmd(\cK)$.
We prove that the stated isomorphism is given by
\begin{eqnarray*}
\Entw^w(\cK)(Y(l),(t,c,\psi)) &\to&
\Cmd(\cK)(I_*(l),\Cmd(Q)((c,\psi),\delta,\varepsilon)),\\
(W,\varrho,\kappa) \1c{\omega} (W',\varrho',\kappa')
&\mapsto &
\mathrm{Cmd}(Q)((W,\varrho),\kappa) \1c{Q(\omega)} 
\mathrm{Cmd}(Q)((W',\varrho'),\kappa').
\end{eqnarray*}
If applying the convention of choosing trivial splittings of identity
2-cells, as described in Theorem \ref{thm:Mnd^i}, then when restricted to the
2-subcategory $\Mnd(\cK)$ of $\Mnd^\iota(\cK)$, $Q$ is equal to $J$. Hence
by \cite[Theorem 2]{Street}, there is a category isomorphism
\begin{eqnarray}\label{eq:2-adj}
&\ \ \cK(l,Q(t))\to \Mnd(\cK)(I(l),t), \ \  
&V \1c \omega V'\ \mapsto\  (vV,v\epsilon V) \1c {v\omega} (vV',v \epsilon
V'); \\
&\ \ \Mnd(\cK)(I(l),t) \to \cK(l,Q(t)), \ \ 
&(W,\varrho) \1c \varphi (W',\varrho') \ \mapsto \ Q(W,\varrho)
\1c{Q(\varphi)} Q(W',\varrho').
\nonumber
\end{eqnarray}
We claim that there is a bijection also between 2-cells $(W,\varrho)\2c \kappa
(c,\psi)(W,\varrho)$ in $\Mnd^\iota(\cK)$, and 2-cells
$Q(W,\varrho)\2c \xi Q(c,\psi)Q(W,\varrho)$ in $\cK$, for any
1-cell $I(l) \1c{(W,\varrho)}t$ in $\Mnd(\cK)$. 
Indeed, for a 2-cell $\kappa$ as described, 
$\xi:=\big( Q(W,\varrho) \2c{Q(\kappa)} Q((c,\psi)(W,\varrho)) \2c \cong
Q(c,\psi)Q(W,\varrho)\big)$ is a 2-cell in $\cK$ as needed. Conversely, for a
2-cell $\xi$ as above, use the chosen splitting $cv \2c \pi v Q(c,\psi) \2c
\iota cv$ of the idempotent 2-cell \eqref{eq:idemp} to construct a 2-cell
$\kappa:=\iota Q(W,\varrho)\ast v \xi:W \2c{} cW$ in $\cK$. It satisfies 
\begin{eqnarray*}
\kappa \ast \varrho
&=& \iota Q(W,\varrho) \ast v \xi \ast v \epsilon Q(W,\varrho)
= \iota Q(W,\varrho) \ast v\epsilon Q(c,\psi)Q(W,\varrho) \ast tv \xi
\nonumber\\ 
&\stackrel{\eqref{eq:Q_1-cell}}{=}& 
\iota Q(W,\varrho) \ast \pi Q(W,\varrho) \ast c\varrho \ast \psi W \ast
t\iota Q(W,\varrho) 
\ast tv \xi  
= c\varrho \ast \psi W \ast t \kappa, 
\end{eqnarray*}
where the last equality follows by 
$\iota f\ast \pi f \ast \psi 
= c\mu  \ast \psi t \ast \eta ct \ast \psi 
\stackrel{\eqref{eq:wentw_m}}{=} \psi \ast \mu c \ast \eta tc 
= \psi$.
Hence $\kappa$ is a 2-cell $(W,\varrho)\2c{} (c,\psi)(W,\varrho)$ in
$\Mnd^\iota(\cK)$, as required. 
 In order to see that this correspondence $\kappa \leftrightarrow \xi$ is
a bijection, note that by \eqref{eq:Q_2-cell},
$vQ(\iota Q(W,\varrho))$ is equal to the composite of $vQ(c,\psi)Q(W,\varrho)
\2c{\iota Q(W,\varrho)} cW$ and the chosen epi 2-cell $cW \2c {}
vQ((c,\psi)(W,\varrho))$. That is,
$Q(\iota Q(W,\varrho))$ is equal to the coherence iso 2-cell
$Q(c,\psi)Q(W,\varrho) \2c \cong Q((c,\psi)(W,\varrho))$. Hence starting
with a 2-cell $\xi$ and iterating both constructions, we re-obtain
$\xi$. 
In the opposite order, applying both constructions to $\kappa$, by
\eqref{eq:Q_2-cell} we get 
$
\iota Q(W,\varrho) \ast \pi Q(W,\varrho) \ast \kappa$.  
This is equal to $\kappa$ by 
\begin{equation}\label{eq:ip_kap}
\iota Q(W,\varrho) \ast \pi Q(W,\varrho) \ast \kappa 
=c \varrho \ast \psi W \ast \eta cW \ast \kappa
\stackrel{\eqref{eq:entw_mod}}{=} \kappa.
\end{equation}
Next we show that $Q(W,\varrho)\2c{Q(\kappa)} 
Q((c,\psi)(W,\varrho)) \2c \cong 
Q(c,\psi)Q(W,\varrho)$ is a
coassociative coaction if and only if $W\2c \kappa cW$ is coassociative, and
it is counital if and only if $\kappa$ is counital.
Compose the coassociativity condition $Q((c,\psi)\kappa) \ast Q(\kappa) =
Q(\delta (W,\varrho)) \ast Q(\kappa)$ horizontally by $v$ on the left and
compose it vertically by the chosen mono 2-cell
$vQ((c,\psi)(c,\psi)(W,\varrho)) \2c \iota ccW$ on the left. 
Applying \eqref{eq:Q_2-cell}, \eqref{eq:entw_mod} and \eqref{eq:ip_kap},
the resulting equivalent condition can be written in the form
$c \kappa \ast \kappa 
= cc \varrho \ast c \psi W \ast \psi c W \ast \eta cc W \ast \delta W\ast
\kappa$. Since 
\begin{equation*}
cc \varrho \ast c \psi W \ast \psi c W \ast \eta cc W \ast \delta W \ast
\kappa \stackrel{\eqref{eq:wentw_d}}{=}
\delta W\ast c\varrho \ast \psi W \ast \eta cW \ast \kappa
\stackrel{\eqref{eq:ip_kap}}{=} \delta W \ast \kappa,
\end{equation*}
this proves that the coaction on $Q(W,\varrho)$ is coassociative if and only
if $\kappa$ is so. 
By \eqref{eq:2-adj}, \eqref{eq:Q_2-cell} and \eqref{eq:ip_kap},
the counitality condition 
$Q(\varepsilon (W,\varrho)) \ast Q(\kappa) =Q(W,\varrho)$ 
is equivalent to $\varepsilon W \ast \kappa=W$.
Thus there is a bijection between the objects of 
$\Cmd(\cK)(I_*(l),\Cmd(Q)((c,\psi),\delta,\varepsilon))$
and the objects of $\Entw^w(\cK)(Y(l),(t,c,\psi))$, as stated.

One can see by similar steps that, for a 2-cell $(W,\varrho) \2c \omega  
(W',\varrho')$ in $\Mnd(\cK)$, $Q(\omega)$ is a morphism $Q(W,\varrho)\to
Q(W',\varrho')$ in $\Cmd(\cK)(I_*(l),\Cmd(Q)((c,\psi),\delta,\varepsilon))$
if and only if  
$\kappa'\ast \omega = c \varrho' \ast \psi W'\ast \eta c W' \ast c \omega \ast
\kappa$. Since
\begin{equation*}
c \varrho' \ast \psi W'\ast \eta c W' \ast c \omega \ast \kappa
=c \varrho' \ast ct \omega \ast \psi W \ast t\kappa \ast \eta W
=c \omega \ast c \varrho \ast \psi W \ast t\kappa \ast \eta W 
\stackrel{\eqref{eq:ip_kap}}{=} c\omega \ast \kappa,
\end{equation*}
we conclude that $Q(\omega)$ is a morphism of $\cK(l,Q(c,\psi))$-coalgebras 
as needed, if and only if $\omega$ is a 1-cell $I_*(l)\1c{} c$ in $\Cmd(\cK)$,
i.e. $\omega$ is a morphism $(W,\varrho,\kappa) \to (W',\varrho',\kappa')$
in $\Entw^w(\cK)(Y(l),(t,c,\psi))$. In view of the isomorphism \eqref{eq:2-adj},
this proves the stated isomorphism
$\Cmd(\cK)(I_*(l),\Cmd(Q)((c,\psi),\delta,\varepsilon))\cong
\Entw^w(\cK)(Y(l),(t,c,\psi))$.

There is a pseudo-natural transformation 
\begin{equation}\label{eq:ps_2-cell}
\Entw^w(Y(-),-) \to
\mathrm{Cmd}(\cK)(\mathrm{Cmd}(Q)AY(-),\mathrm{Cmd}(Q)A(-)), 
\end{equation}
with 1-cell parts the functors induced by the pseudo-functor
$\mathrm{Cmd}(Q)A$ and 2-cell parts provided by its pseudo-naturality
isomorphisms. 
Recall that $AY$ differs from $\mathrm{Cmd}(I)I_*$ by the inclusion 2-functor 
$\mathrm{Cmd}(\mathrm{Mnd}(\cK)) \hookrightarrow
\mathrm{Cmd}(\mathrm{Mnd}^\iota(\cK))$.
Since applying $Q:\mathrm{Mnd}^\iota(\cK)\to \cK$ after $\cK \1c I
\mathrm{Mnd}(\cK)\hookrightarrow \mathrm{Mnd}^\iota(\cK)$ we obtain the 
identity functor $JI=\cK$, it follows that $\mathrm{Cmd}(Q)AY(-)=I_*$ as
pseudo-functors. Thus \eqref{eq:ps_2-cell} is, in fact, a pseudo-natural
transformation $\Entw^w(Y(-),-) \to
\mathrm{Cmd}(\cK)(I_*(-),\mathrm{Cmd}(Q)A(-))$. 
Since we already proved that its 1-cells are isomorphisms, it is a
pseudo-natural isomorphism, as stated.  
\end{proof}

\begin{theorem}\label{thm:EM_iso}
Let $\cK$ be a 2-category which admits Eilenberg-Moore constructions for both
monads and comonads and in which idempotent 2-cells split. The following
diagram of pseudo-functors is commutative, up to a pseudo-natural equivalence.  
$$
\xymatrix{
\Entw^w(\cK)\ar[rrr]^-{A}\ar[dd]_-{B}&&&
\Cmd(\Mnd^\iota(\cK))\ar[d]^-{\Cmd(Q)}\\
&&&\Cmd(\cK)\ar[d]^-{J_*}\\
\Mnd(\Cmd^\pi(\cK))\ar[rr]^-{\Mnd(Q_*)}&&
\Mnd(\cK)\ar[r]^-{J}&
\cK \ .}
$$
In particular, for any weak entwining
structure $(t,c,\psi)$ in $\cK$, the monad
$\mathrm{Mnd}(Q_*)(c\1c{(t,\psi)}c,\mu,\eta)$ and the comonad 
$\mathrm{Cmd}(Q)(t\1c{(c,\psi)}t,\delta,\varepsilon)$ in $\cK$ possess
equivalent Eilenberg-Moore objects. 
\end{theorem}

\begin{proof}
The proof consists of showing that both $J_* \Cmd(Q)A$ and $J\Mnd(Q_*)B$ are
right biadjoints of the 2-functor $Y$ in Corollary \ref{cor:functors}(1).
Then the claim follows by uniqueness of a biadjoint up to a pseudo-natural
equivalence.  

On one hand, there is a sequence of pseudo-natural isomorphisms 
$$
\cK(-,J_* \Cmd(Q)A(-))
\cong \Cmd(\cK)(I_*(-),\Cmd(Q)A(-))
\cong \Entw^w(\cK)(Y(-),-),
$$
where the second isomorphism follows by Proposition \ref{prop:comodules}.

On the other hand, applying Proposition \ref{prop:comodules} to the 2-category
$\cK_*$ (in the third step) and using Corollary \ref{cor:functors}(2) (in the
last step), we obtain a sequence of pseudo-natural isomorphisms 
\begin{eqnarray*}
\cK(-,J \Mnd(Q_*)B(-))
&\cong& \Mnd(\cK)(I(-),\Mnd(Q_*)B(-))\\
&\cong& \Cmd(\cK_*)_*(I(-),\Mnd(Q_*)B(-))\\
&\cong& \Entw^w(\cK_*)_*(\Phi Y(-),\Phi(-))
\cong \Entw^w(\cK)(Y(-),-).
\end{eqnarray*}
\end{proof}

\begin{example}
Consider the 2-subcategory $\cK$ of $\mathrm{CAT}$, whose 1-cells are functors
induced by bimodules. Explicitly, 0-cells be module categories $M_R$ for
algebras $R$ over a fixed commutative ring $k$. The 1-cells $M_R \1c{} M_{R'}$
be $R$-$R'$ bimodules $V$, i.e. functors $(-)\ot_R V:M_R\to M_{R'}$. The
2-cells $V\2c{} W$ be $R$-$R'$ bimodule maps $\omega:V\to W$, i.e. natural
transformations $(-)\ot_R V \2c{(-)\ot_R \omega} (-)\ot_R W$.  

A weak entwining structure in $\cK$ is then a triple $(t:=(-)\ot_R T,
c:=(-)\ot_R C,\psi:=(-)\ot_R \Psi)$, where $R$ is a $k$-algebra, $T$ is an
$R$-ring (i.e. a monad $R\1c T R$ in $\mathrm{BIM}_k$), $C$ is an $R$-coring
(i.e. a  comonad $R\1c C R$ in $\mathrm{BIM}_k$), and $\Psi:C\ot_R T \to T\ot_R
C$ is an $R$-bimodule map such that the equalities
(\ref{eq:wentw_m}-\ref{eq:wentw_e}) hold true. 

In this particular 2-category $\cK$, the idempotent 2-cell
\eqref{eq:idemp} is given by an idempotent map. Taking its obvious splitting
through its range, the associated pseudo-functor $Q: \mathrm{Mnd}^\iota(\cK)
\to \cK$ in Theorem \ref{thm:Mnd^i} becomes a 2-functor. Hence the
isomorphisms in Proposition \ref{prop:comodules} become 2-natural, so that the
equivalent Eilenberg-Moore objects in Theorem \ref{thm:EM_iso} become
isomorphic.  

Under the minor restriction that $R=k$, the monad $\Mnd(Q_*)B((-)\ot_R
T,(-)\ot_R C,$ $(-)\ot_R \Psi)$ and the comonad $\Cmd(Q)A((-)\ot_R T,(-)\ot_R
C,(-)\ot_R \Psi)$ were described in \cite[Section 2]{CaeDeG}. It was shown in
\cite[Proposition 2.3]{Brz:coring} that their Eilenberg-Moore categories are
isomorphic to the category of so-called weak entwining structures. Using the
constructions in the current paper, this category of weak entwining structures
is nothing but $\Entw^w(\cK)(Y(k),((-)\ot_R T,(-)\ot_R C,(-)\ot_R \Psi))$.  
\end{example}

\section*{Acknowledgement}
The author's work is supported by the Hungarian Scientific Research Fund OTKA
F67910. 
She is grateful to the referee of \cite{weak_EM}, whose comments were
highly valuable also for this work.

\end{document}